\documentclass[reqno]{amsart}

\usepackage{amsmath, amsfonts, amsthm, amssymb, color,  graphicx, mathrsfs, cite}
\usepackage{comment,stmaryrd}

\theoremstyle{plain}
\newtheorem{theorem}{Theorem}[section]
\newtheorem{lemma}{Lemma}[section]
\newtheorem{proposition}{Proposition}[section]
\newtheorem{corollary}{Corollary}[section]

\theoremstyle{definition}
\newtheorem{definition}{Definition}[section]

\theoremstyle{remark}

\newtheorem{example}{Example}[section]

\numberwithin{equation}{section}
\allowdisplaybreaks

\usepackage{ifpdf}
\ifpdf \usepackage[colorlinks=true, citecolor=blue, linkcolor=blue, urlcolor=blue]{hyperref} \fi

\def\thrm{\begin{theorem}}
\def\thrml#1{\begin{theorem}\label{#1}}
\def\ethrm{\end{theorem}}
\def\lmm{\begin{lemma}}
\def\lmml#1{\begin{lemma}\label{#1}}
\def\elmm{\end{lemma}}
\def\dfntn{\begin{definition}}
\def\dfntnl#1{\begin{definition}\label{#1}}
\def\edfntn{\end{definition}}
\def\crllr{\begin{corollary}}
\def\crllrl#1{\begin{corollary}\label{#1}}
\def\ecrllr{\end{corollary}}
\def\xmpl{\begin{example}}
\def\xmpll#1{\begin{example}\label{#1}}
\def\exmpl{\end{example}}
\def\nmrt{\begin{enumerate}}
\def\enmrt{\end{enumerate}}
\def\qtn{\begin{equation}}
\def\qtnl#1{\begin{equation}\label{#1}}
\def\eqtn{\end{equation}}
\def\prpstn{\begin{proposition}}
\def\prpstnl#1{\begin{proposition}\label{#1}}
\def\eprpstn{\end{proposition}}

\def\tm#1{\item[{\rm (#1)}]}
\def\proof{{\bf Proof}.\ }
\def\eprf{\hfill$\square$}

\def\cX{\mathcal {X}}

\begin{document}

\title{The endomorphism rings of permutation modules of $\frac{3}{2}$-transitive permutation groups}
\maketitle

\begin{center}
{\author Jiawei He$^a$ and Xiaogang Li$^{b,*}$}
\end{center}

\vskip 3mm

\begin{abstract}
Recent classification of $\frac{3}{2}$-transitive permutation groups leaves us with six families of groups which are $2$-transitive, or Frobenius, or one-dimensional affine, or the affine solvable subgroups of $ \mathrm{AGL}(2, q)$, or special projective linear group $\mathrm{PSL}(2, q)$, or $\mathrm{P\Gamma L}(2, q)$, where $q=2^p $ with $p$ prime. According to a case by case analysis, we prove that the endomorphism ring of the natural permutation module for a $\frac{3}{2}$-transitive permutation group is a symmetric algebra.

{ Keywords: Permutation module; Endomorphism ring; Symmetric algebra.}

\end{abstract}


\date{}

\maketitle


\renewcommand{\thefootnote}{\empty}
\footnotetext{School of Mathematics and Information Science, Nanchang Hangkong University, Nanchang, China}
\footnotetext{$^*$corresponding author, a School of Mathematical Sciences, Capital Normal University, Beijing, 100048, P.R.China}
\footnotetext{ Email:  hjwywh@mails.ccnu.edu.cn }

\section{Introduction}\label{in}

The concept of $\frac{3}{2}$-transitive permutation group goes back to research of Wielandt [10], related with constructing of many combination objects; the defining property of such a group is a non-regular transitive group such that all the nontrivial orbits of a point stabilizer have equal size. More steps towards the classification of the primitive $\frac{3}{2}$-transitive groups were taken in [1] and [4]. It was proved that primitive $\frac{3}{2}$-transitive groups are either affine or almost simple, and the almost simple examples were determined. Further more, if we do not require primitivity of a $\frac{3}{2}$-transitive permutation groups,  a recent classification of $\frac{3}{2}$-transitive permutation group can be found in [6].

\medskip
Given a field $k$ and a natural number $n$, if we assume that $G$ is a permutation subgroup of the symmetric group on $\Omega$ with $n$ letters, then $G$ can be naturally viewed as a subgroup of permutation matrices in the general linear group $\operatorname{GL}(n,k)$. It is well known that the {\it centralizer ring} which is defined to be the centralizer of $G$ in the full matrix ring $M_n(k)$ is isomorphic to the endomorphism ring of the permutation $kG$-module $k\Omega$. An important research on the centralizer ring theory of Wielandt [10] is to study the centralizer ring corresponding to the permutation group of $kp$ degree for a prime $p$. However, when $G\leq \rm Sym (\Omega)$ is a $\frac{3}{2}$-transitive permutation group, some classical results on permutation groups can be applied to the centralizer ring (or the endomorphism ring). In this paper, we study the endomorphism ring of the natural permutation modules of a $\frac{3}{2}$-transitive groups. The following is our main theorem:
\medskip

\thrml{1451b} Let $k$ be a field and $G$ be a $\frac{3}{2}$-transitive permutation groups on a set $\Omega$, then the endomorphism ring of the natural permutation module $k\Omega$ for $kG$ is a symmetric algebra.
\ethrm
\medskip

Theorem 1.1 is proved in Section 3. See Subsection 2.2 for related definitions.

\medskip

Throughout this paper, groups and modules are always finite groups and right modules. For a finite set $X$, a group $G$ and a subgroup $H\leq G$, we always use $|X|$, $|G|$ and $|G:H|$ to denote the cardinality of $X$, the order of $G$ and the index of $H$ in $G$ respectively. Let $k$ be a field, then we always use $k_H$ and $k^G_H$ to denote the trivial $kH$-module and the induced module $k_H\otimes_{kH} kG$ respectively. For two $kG$-modules $M$ and $N$, we use $M|N$ to mean that $M$ is isomorphic to a direct summand of $N$ as a $kG$-module. If moreover $G$ is a permutation group on some set $\Omega$, then the $k$-linear span of elements of $\Omega$ which we denote by $k\Omega$ is a natural $kG$-module via the action $$(\sum_{\alpha\in \Omega} a_{\alpha}\alpha) \cdot g:=\sum_{\alpha\in \Omega} a_{\alpha}(\alpha\cdot g), \forall\ g\in G.$$
If we denote $\underline{H}=\sum_{h\in H} h$ and $\underline{G}=\sum_{g\in G} g$, then it can be easily checked that $k\Omega\cong \underline{H}kG$ as $kG$-module. If ${\rm char} ~k\nmid |G:H|$, then $$\underline{H}kG=\underline{G}kG\oplus(\underline{H}-|G:H|^{-1}\underline{G})kG,$$ which implies that $k_G\cong \underline{G}kG$ is isomorphic to a direct summand of $k\Omega$ in this situation.

\medskip

$\mathbf{Notation.}$

Throughout this paper, $\Omega$ is always a finite set. $A_{|\Omega|}$ and $S_{|\Omega|}$ denote the alternating group and symmetric group on $\Omega$, respectively.

The diagonal of the Cartesian product $\Omega \times \Omega$ is always denoted by $1_{\Omega} ;$ for $\Delta\subseteq\Omega$, we set $1_{\Delta}:=\{(\alpha,\alpha):\alpha\in \Delta\}$.

For a subset $r \subseteq \Omega \times \Omega$ and any $\gamma \in \Omega$, we set $r^{*}=\{(\beta, \alpha):(\alpha, \beta) \in r\}$, $\gamma r=\{\delta \in \Omega:(\gamma, \delta) \in r\}$  and $\Omega(r):=\{\alpha\in \Omega: \alpha s\neq \emptyset \}$.

\medskip

\section{Preliminary}\label{Pr}

\subsection{Symmetric Algebras}

\dfntnl{1055a} Let $k$ be a field. A $k$-algebra $A$ is called {\it self-injective} if every finitely generated projective $A$-module is injective.
\edfntn

 We shall develop some of the properties of these algebras here. A special class of self-injective algebras plays an important role as
group algebras are of this kind.

\dfntnl{1055a} Let $k$ be a field.  A $k$-algebra $A$ is called $symmetric$ if $\operatorname{Hom}_{k}(A, k) \cong A$ as $A$-$A$-bimodules, or equivalently, there exists a non-degenerate symmetric associative bilinear form $f: A \times A \rightarrow k$.
\edfntn

We refer the readers to [9, Chapter 2, Section 8] for properties of symmetric algebras.

\subsection{Coherent Configurations}

\dfntnl{1055a}  A pair $\cX=(\Omega, S)$ where $S$ is a partition of $\Omega\times \Omega$ is called {\it coherent configuration} if
\nmrt
  \tm{1} $1_{\Omega} \in S^{\cup}$, where elements of $S^{\cup}$ are unions of several elements in $S$,
  \tm{2}  $S^{*} = S$, where $S^{*}:=\{r^{*}:r\in S\}$,
  \tm{3}  For any $r, s, t\in S$, the number $C^t_{rs} := |\alpha r\cap \beta s^{*}|$ does not depend on the choice of $(\alpha, \beta)\in t$, here $|\alpha r\cap \beta s^{*}|$ is the number of $\gamma\in \Omega$ such that $(\alpha,\gamma)\in r, ~(\gamma, \beta)\in s$.
  \enmrt
\edfntn

Each element in the set $S$ is called a {\it primitive set}. For any $\delta \in \Omega ({s})$, the positive integer $|\delta s|$ equals to the intersection number $c_{s s^{*}}^{1_{\Omega ({s})}}$, hence does not depend on the choice of $\delta$ and we call it the {\it valency} of $s$. By [2], the following equalities hold:
$$
|t| C_{r s}^{t^{*}}=|r| C_{s t}^{r^{*}}=|s| C_{t r}^{s^{*}}, \quad r, s, t \in S.
$$
\medskip

For a permutation group $G\leq \rm {Sym}(\Omega)$, by [2], we can obtain a coherent configuration associated to the group $G$. All 2-orbits under the action of $G$ on $\Omega\times \Omega$ are the primitive sets of this coherent configuration. For any $s\in {\rm {Orb}}(G, \Omega^2)$, it is said to be {\it reflexive} if $(\alpha,\alpha)\in s$ for some $\alpha\in \Omega$.
\medskip

\subsection{Schur rings}

\dfntnl{1055a}
Given a field $k$ and a finite group $G$. A subring $\mathcal{A}$ of the ring $kG$ is called the $Schur$ $ring$  over $G$ if it has a linear $k$-base consisting of elements $V$, where $V$ runs over a family $\mathcal{P} $ of pairwise disjoint nonempty subsets of $G$ such that
\nmrt
  \tm{1} $\{e\} \in \mathcal{P}$,
\tm{2} $\bigcup\limits_{V\in \mathcal{P}}V=G$,
\tm{3} for any~$V\in \mathcal{P}$, $V^{-1}=\{v^{-1} : v\in V\}\in \mathcal{P}$.
\enmrt
\edfntn

The elements of $\mathcal{P}$ are called the $basic$ $sets$ of $\mathcal{A}$. For any $V\subseteq  G$,  denote $\underline{V}=\sum_{v\in V} v$.
 The following lemma is well-known, which is [8, Theorem 1.6].

\lmml{1440a} Let $G$ be a finite group and $k$ a field. Suppose $T \leq \operatorname{Aut}(G)$. Then the linear space
$$
\mathcal{A}=\mathrm{span}_k\left\{\underline{V} \in k G: V\in \operatorname{Orb}(T, G) \right\}
$$
is an Schur ring over $G .$
\elmm
\section{Proofs}

Note that a permutation group $G\leq \rm {Sym} (\Omega)$ is said to be {\it half-transitive} if all orbits of $G$ are of the same cardinality. Moreover, if $G$ acts transitively on $\Omega$, then the cardinalities of the orbits of $G_\alpha$ for any given point $\alpha\in\Omega$ are called the $subdegrees$ of $G$ and $G$ is called $primitive$ if a point stabilizer in $G$ is a maximal subgroup. Before we proceed on to prove the main theorem, we first prove the following lemma.
\vskip 3mm

\lmml{1440a} Let $k$ be a field and $G$ a half-transitive permutation group on $\Omega$. Suppose that the orbit sizes for the point stabilizer $G_{\alpha}$ on $\Omega$ are invertible in $k$ for all $\alpha\in \Omega$, then ${\rm End}_{kG}(k\Omega)$ is a symmetric algebra.

\elmm

\proof
Note that $M_n(k)^G={\rm End}_{kG}(k\Omega)$ when regarding $G$ as permutation matrices, thus we can use matrices to describe ${\rm End}_{kG}(k\Omega)$. We take $\Omega$ as the set $\{1,2,\cdots,n\}$ for convenience. Now we denote the orbits under the action of $G$ on $\Omega\times \Omega$ by
$$S=\{s_1,\cdots, s_t, s_{t+1}, \cdots, s_m\},$$
where $\{s_1,\cdots,s_t\}$ are precisely the reflexive orbits. For any orbit $s_l\in S$, we define a matrix as follows:
\[
(A_{s_l})_{ij}=\left\{
\begin{array}{rcl}
1   &    &{(i,j)\in s_l}\\
0   &    &{(i,j)\notin s_l}\\
\end{array}
\right.
\]
called the adjacency matrix of $s_l$. It is known from [2, Theorem 2.3.5] that $$M_n(k)^G=\mathop  \oplus \limits_{i = 1}^m kA_{s_i}$$ as vector spaces. We fix $\mathcal{B}=\{A_{s_i},1\leq i\leq m\}$ as a basis of $M_n(k)^G$. For $s_l\in S$, we conclude from [2, Section 1.1] that $\Omega(s_l)$ is an orbit of the action of $G$ on $\Omega$.
We now define a map $$\varphi: M_n(k)^G\times M_n(k)^G\rightarrow k,~(A_{s_l}, A_{s_{l'}})\mapsto C^{1_{\Omega(s_l)}}_{s_l s_{l'}}1_k$$ for all $A_{s_l},A_{s_{l'}}\in \mathcal{B}$, where $1_k$ is the identity of $k$. From [2, Chapter 1], $\varphi$ is clearly bilinear, while $\varphi(A_{s_l}, A_{s_{l'}}A_{s_{l''}})$ and $\varphi(A_{s_l}A_{s_{l'}}, A_{s_{l''}})$ are just the coefficient of $A_{1_{\Omega(s_l)}}$ in the expression of $A_{s_l} A_{s_{l'}}A_{s_{l''}}$ with respect to the basis $\mathcal{B}$, they must be equal and therefore $\varphi$ is associative. Then we only need to show that $\varphi$ is symmetric and non-degenerate. To show it is symmetric, it suffices to show that for all $s_l,s_{l'}\in S$, the following equality holds: $$C^{1_{\Omega(s_l)}}_{s_l s_{l'}}1_k=C^{1_{\Omega(s_{l'})}}_{s_{l'} s_l}1_k.$$ Since $G$ is half-transitive on $\Omega$, ie., all orbits of $G$ on $\Omega$ have the same size, but from [2, Chapter 1], we have $$C^{1_{\Omega(s_{l'})}}_{s_{l'} s_l}|\Omega(s_{l'})|1_k=C^{1_{\Omega(s_l)}}_{s_l s_{l'}}|\Omega(s_{l})|1_k,$$ thus we obtain the required equation.

\vskip 3mm

Now it remains to show that $\varphi$ is non-degenerate, denote
$$\mathrm{Rad}\varphi:=\{x\in M_n(k)^G:\varphi(x,y)=0~ \text{for~all} ~y\in M_n(k)^G\},$$
we have to show that Rad$\varphi=0$. Now choose arbitrary $x\in$ Rad$\varphi$, we can write $x=\sum^m_{i=1} a_{s_i} A_{s_i}$, where $a_{s_i}\in k$ for all $1\leq i\leq m$. For any $s_l\in S$, it is clear that
$$C^{1_{\Omega(s_l)}}_{s_l s_{l'}}1_k\neq 0 \Leftrightarrow s_l=s^{*}_{l'}.$$
Choose an element $(a,b)\in s_l$ and denote $G_a$ the stabilizer for $a$ in $G$, then we deduce that $$C^{1_{\Omega(s_l)}}_{s_l s_{l'}}1_k=|b^{G_a}|1_k$$
which is an orbit size under the action of $G_a$ on $\Omega$. Therefore we have the following equations
$$\varphi(x, s^{*}_l)=a_{s_l} C^{1_{\Omega(s_l)}}_{s_l s^{*}_l}=0,$$
which implies that $a_{s_l}=0$ . Hence we derive that $a_{s_i}=0$ when $s^{*}_i$ runs over the orbits of $G$ on $\Omega\times \Omega$, therefore $a=0$, as required.
\eprf
\medskip

\xmpl
Given a field $k$ with characteristic $p$ and a natural number $n\geq 3$, we can define a map $\varphi$ from $S_n$ to $A_{n+2}$ by the following rule:
\[
\varphi(a):=\left\{
\begin{array}{ll}
a, &  \text{if}~{\rm sgn}(a)=1 \\
a(n+1,n+2), & \text{otherwise},
\end{array}
\right.
\]
here sgn is the Sign Function. Then $\varphi$ is an injection. We now view $S_n$ as a subgroup of $A_{n+2}$ in this way, and let $A_{n+2}$ act on the the set $\Omega:=\{S_{n}x:x\in A_{n+2}\}$ of all right cosets of $S_n$, then this action is faithful and primitive. Moreover, $A_{n+2}$ has degree $\frac{(n+1)(n+2)}{2}$ and three subdegrees $$\{1,2n,\frac{n(n-1)}{2}\}.$$ We claim that if $p\geq 5$, then ${\rm End}_{kA_{n+2}}(k\Omega)$ is always a symmetric algebra. To see it, note that $p$ only divides one of $n-1, n, n+1, n+2$, we may assume $p\mid n(n-1)$ by Lemma 3.1. But now $$k^{A_{n+2}}_{S_n}:=k\otimes_{k{S_{n}}}k{A_{n+2}}$$ has the trivial module $k$ as a summand. Note that the rank of $A_{n+2}$ is $3$ means that ${\rm End}_{kA_{n+2}}(k\Omega)$ is 3-dimensional. If we write $k\Omega=k\oplus N$ for some $kG$-module $N$, then ${\rm End}_{kA_{n+2}}(k\Omega)\cong k \times {\rm End}_{kA_{n+2}}(N)$. This implies that ${\rm End}_{kA_{n+2}}(N)$ must be 2-dimensional. Thus ${\rm End}_{kA_{n+2}}(N)$ is a symmetric algebra and so is ${\rm End}_{kA_{n+2}}(k\Omega)$.
\exmpl
\vskip 3mm

{\bf Proof of Theorem \ref{1451b}.}\, We divide our proof into four parts according to the four situations in the classification of $\frac{3}{2}$-transitive permutation groups.
\vskip 3mm

\vskip 3mm

(1) {\bf $G$ is 2-transitive.}

\medskip

Denote $|\Omega|=n$, then $${\rm End}_{kG}(k\Omega)\cong k[X]/(X^2)~ \text{or} ~{\rm End}_{kG}(k\Omega)\cong k\times k$$ depending whether ${\rm char}~k$ divides $n$ or not, ${\rm End}_{kG}(k\Omega)$ is always a symmetric algebra in both cases.

\vskip 3mm

(2) {\bf $G$ is a Frobenius group.}

\medskip

Denote $|\Omega|=n$ and ${\rm char}~k=p$, let $\alpha\in \Omega$ and $H=G_{\alpha}$ be the point stabilizer of $\alpha$ in $G$. We may assume that 
$p\mid |G|$, otherwise ${\rm End}_{kG}(k\Omega)$ would be semisimple by Maschke's theorem. That $G$ is a Frobenius group implies that $|H|$ and $n$ are coprime, thus either $p\mid |H|$ or $p\mid n$.
\medskip

We first suppose that $p\mid |H|$. Let $P$ be a Sylow $p$-subgroup of $H$, then $N_G(P)\leq H$. Now $k\Omega\cong k^G_H$ as $kG$-module and $P$ is a vertex of $k_H$, hence the Scott module $k_G$ is the Green correspondence of $k_H$ and the vertice of any indecomposable summands of $k^G_H$ except $k_G$ is of the form $P\cap P^g$ for some $g\in G-H$. But that $G$ is a Frobenius group yields that $H$ has trivial intersection property, therefore the vertices of all indecomposable summands of $k^G_H$ except $k_G$ are the identity and hence those summands are projective. Whence we can write $k\Omega=k_G\oplus M$ for some projective $kG$-module and $${\rm End}_{kG}(k\Omega)\cong k\times {\rm End}_{kG}(P).$$ It is well-known that $kG$ is a symmetric algebra and thus so is ${\rm End}_{kG}(P)$.

\medskip
If $p\mid n$, then $|H|$ is invertible in $k$ and we set $$e:=\frac{1}{|H|}\sum_{h\in H}h.$$ Then direct computation shows that $k\Omega\cong ekG$ as $kG$-module, and then $${\rm End}_{kG}(k\Omega)\cong ekGe$$ is a symmetric algebra.

\medskip
In fact, we have the following more general result using the concept of Schur rings:

\thrml{14512b}  Let $k$ be a field and $G$ a transitive permutation group on $\Omega$ with a regular normal subgroup $T$. Let $H=G_{\alpha}$ be a point stabilizer for $\alpha\in \Omega$ and suppose that all subdegree of $G$ are coprime to $|T|$ and $p\nmid |H_{\beta}|$ for some $\beta\in\Omega-\{\alpha\}$, then ${\rm End}_{kG}(k\Omega)$ is always a symmetric algebra.

\ethrm
\medskip
\proof
By Lemma 3.1, it suffices to assume that ${\rm char}\ k$ divides some subdegrees of $G$. Fixed some $\alpha\in \Omega$, for any $t\in T, g\in G$, we define $t^g$ to be the unique element in $T$ such that $\alpha\cdot t^g=\alpha\cdot (tg)$. It can be directly checked that this define a group action of $G$ on $T$ which is permutation isomorphic to the group action of $G$ on $\Omega$. As a consequence of normality of $T$, we have $x^h=h^{-1}xh$ for all $x\in T, h\in H$. Denote $O_1,\cdots,O_s$ by the orbits under the action of $H$ on $T$. Then by Lemma 2.1,
 $$\mathcal{A}:={\rm span}_k\{\underline{O_i}: 1\leq i\leq s\}$$
  a Schur ring $\mathcal{A}$ over $T$.
Invoking Schur's method, we know that ${\rm End}_{kG}(k\Omega)$ is isomorphic to $\mathcal{A}$. We have ${\rm char}\ k\nmid |T|$ by assumption, which yields that $kT$ is a semisimple algebra by Maschke's theorem. Let $t\in T$  be the element such that $\alpha\cdot t=\beta$ and denote $O_t$ by the orbit containing $t$, then
$$\sum_{h\in H} t^h=|H_{\beta}|\underline{O_t}$$
is nonzero in $kT$, where $H_\beta$ is the stabilizer of $\beta$ in $H$. Note that the conjugation of $H$ on $T$ can be extended to an automorphism of the group algebra $kT$ just via the natural way. Invoking [7, Corollary 5.4], we deduce that $$(kT)^H:=\{x\in kT : x^h=x~\text{for~all}~h\in H\}$$ is also a semisimple algebra. But $(kT)^H$ is exactly this Schur ring $\mathcal{A}$, we are done.
\eprf

\vskip 3mm

(3) {\bf $G$ is an affine primitive permutation group}

\medskip

That is, $G= {\rm T}(V)H\leq \rm AGL(V)$ for some vector space $V$, where $H \leq \rm GL(V)$ is a $\frac{1}{2}$-transitive irreducible linear group, ${\rm T}(V )$ is the group of translations and $V$ is of dimension $d$ over a finite field $\mathbb{F}_p$ for a prime $p$. We denote $T$ by the normal regular subgroup of $G$. This case will follow from the following easy but useful lemma.
\vskip 3mm

 \lmml{1951b} Let $k$ be a field and $G$ a transitive permutation group on some set $\Omega$ with an abelian regular subgroup $T$. If any subdegree of $G$ is coprime to the order of $T$, then ${\rm End}_{kG}(k\Omega)$ is always a symmetric algebra.
\elmm

\vskip 3mm

\proof
By Lemma 3.1, it suffices to consider the situation when the characteristic of $k$ divides some nontrivial subdegree of $G$. We can similarly define an action of $G$ on $T$ as that in the proof of Theorem 3.1. Denote $H$ for some point stabilizer and $O_1,\cdots,O_s$ the orbits under the action of $H$ on $T$. Therefore we know that ${\rm End}_{kG}(k\Omega)$ is isomorphic to a Schur ring $\mathcal{A}$ over $T$, where $$\mathcal{A}:={\rm span}_k\{\underline{O_i}: 1\leq i\leq s\}.$$ Since ${\rm char}\ k$ does not divide the order of $T$, we deduce that $kT$ is a semisimple algebra thanks to Maschke's theorem. As $kT$ is abelian, it has no nonzero nilpotent elements. Hence $\mathcal{A}$ as a subalgebra of $kT$ also contains no nonzero nilpotent elements, forcing $\mathcal{A}$ to be semisimple, we are done.
\eprf

\medskip
Note that if we omit the condition that any subdegree of $G$ is coprime to the order of $T$ in Lemma 3.2, the conclusion may be false, we will give an example in the following.

\medskip

\xmpl
Let $k=\mathbb{F}_3$, we use $X$ to denote the matrix in $M_2(k)$ whose all but the (2,1)-entry are $1$ while the (2,1)-entry is $0$. Let $V=k^2$ and $G\leq \rm AGL(2,3)$ be a linear group such that $G=T\rtimes H$, where $T=\mathbb{Z}_3\times \mathbb{Z}_3$ is the subgroup acting by translation on $V$ and $H$ is the subgroup of $\rm GL(2,3)$ generated by $X$.  We can regard $G$ as a permutation group on $V$ in the natural way and the subdegrees for $G$ are $\{1,1,1,3,3\}$. Invoking Green's indecomposability theorem, we know that $kV\cong k^G_H$ is an indecomposable $kG$-module, therefore the endomorphism ring $E:={\rm End}_{kG}(kV)$ is a local $k$-algebra. Choose suitable generators $a,b$ for $T$. Using Schur ring, we know that $E$ is isomorphic to a $5$-dimensional subalgebra of $kT$ with basis $$\{1,a,a^2,b+ab+a^2b,b^2+ab^2+a^2b^2\}.$$ Now we identify $E$ with that subalgebra. The radical of an $E$-module $R$ which is defined to be the intersection of all maximal submodules of $R$, is denoted by $J(R)$. Note that $J(R)=RJ$ where $J$ is the Jacobson radical of $E$ in our case.
Computations show that $$J(E),~ J^2(E):=J(J(E)),~ J^3(E):=J(J^2(E))$$ are of dimensions $4,2$ and $0$, respectively. In particular, $J^2(E)$ is a semisimple submodule of the regular module $E$. We claim that $E$ is not a symmetric algebra. Suppose for the contrary, then
$\rm soc(E)\cong \rm top(E)$ have to be 1-dimensional, contradicting that $J^2(E)$ is a semisimple $E$-submodule of $E$ whose dimension is $2$.
\exmpl

\medskip
\vskip 3mm

(4) {\bf $G$ is almost simple}

\medskip

For this case, we need the following two lemmas.
\vskip 3mm

\lmml{1951b} Let $k$ be a field, if two permutation groups $G_1\leq G_2$ on $\Omega$ have the same rank, then ${\rm End}_{kG_1}(k\Omega)\simeq {\rm End}_{kG_2}(k\Omega)$.
\elmm

\vskip 3mm

\proof
First, we will show that for a given field $k$ and any permutation group $G$ on some set $\Omega$ with $|\Omega|=n$, the endomorphism algebra of $k\Omega$ is isomorphic to the adjacent algebra of the so-called coherent configuration. Since the natural permutation module $k\Omega$ is in some sense equivalent to the permutation representation $$\rho: kG\rightarrow M_n(k),$$ and to compute the endomorphism ring of $k\Omega$ is to calculate the centralizer ring for permutation matrices corresponding to $G$ in $M_n(k)$, this is $${\rm End}_{kG}(k\Omega)= M_n(k)^G.$$ Now by hypothesis, we have $\rho(G_1)\leq \rho(G_2)$, which implies that $$ {\rm End}_{kG_2}(k\Omega)=M_n(k)^{G_2}\subseteq M_n(k)^{G_1}={\rm End}_{kG_1}(k\Omega).$$ By comparing their rank, we get the equality, as required.
\eprf
\medskip
\lmml{1951b}Let $G$ be a $\frac{3}{2}$-transitive permutation group on $\Omega$ with a point stabilizer $H$, then for any prime $p$ dividing the nontrivial subdegree of $G$, the normalizer of a Sylow $p$-subgroup $P$ of $H$ in $G$ is contained in $H$.
\elmm

\proof
Write $n$ for $|\Omega|$. All the nontrivial subdegrees of $G$ are the same since $G$ is $\frac{3}{2}$-transitive. If $N_G(P)\nleqslant H$, then $$P\leq H\cap H^x$$ for some $x\in N_G(P)-H$ and hence $p\nmid |H:H\cap H^x|$. But $|H:H\cap H^x|$ equals to some subdegree, a contradiction, as required.
\eprf

\medskip

Now we come back to our analysis. In this case, then by [6, Corollary 2.3], we can obtain the following two cases:
 \medskip

(a) $n=21$, $G=A_{7}$ or $S_{7}$ acting on the set of pairs in $\{1, \cdots, 7\}$,

\medskip
(b) $n=\frac{1}{2} q(q-1)$ where $q=2^{m} \geq 8$, and either $G=\operatorname{PSL}(2, q)$, or $G=$ $\mathrm{PGL}(2, \mathrm{q})$ with $m$ prime.
\medskip

For case (a), it is well-known that the ranks of both groups are $3$ and the subdegrees for them are $\{1,10,10\}$.
We may assume that ${\rm char}~k=2$ or $5$ by lemma 3.1. Under the assumption, $$k^G_H=k_G\oplus M$$ for some $kG$-module $M$ and $${\rm End}_{kG}(k\Omega)\simeq k\times {\rm End}_{kG}(M),$$ which implies that ${\rm End}_{kG}(M)$ is 2-dimensional and hence must be a symmetric algebra, thus the same is true for ${\rm End}_{kG}(k\Omega)$.

\medskip
For case (b), it suffices to analyse $G={\rm PSL}(2,q)$ by Lemma 3.3 while the stabilizer $H$ is equal to $D_{2(q+1)}$. Since it suffices to consider the case when ${\rm char}~k$ divides the subdegree $q+1$ of $G$, in particular, ${\rm char}~k$ is an odd prime $p$. Now $$|H\cap H^x|=2$$ for any $x\in G-H$ implies that any Sylow $p$-subgroup $P$ of $H$ has the property that $P\cap P^x=1$ for all $x\in G-H$. Since Lemma 3.4 guarantees that Green's correspondence theorem can be applied, we are done by the same analysis as in situation 2.
\eprf

\medskip

$\mathbf{Remark}$:
Note that by O'Nan-Scott theorem, a primitive permutation group $G$ on $\Omega$ with degree a prime power is either almost simple, or of affine type or of product type. If we suppose that $G$ is almost simple, then by [5, Corollary 2], we know that ${\rm soc}(G)$ which is defined to be the subgroup generated by all minimal normal subgroups of $G$, must act $2$-transitive on $\Omega$ except for the case ${\rm soc}(G)\cong \rm PSU(4,2)$ and $|\Omega|=27$. For the exceptional case, the subdegree for ${\rm soc}(G)$ are $\{1,10,16\}$, but it is not hard to check that ${\rm End}_{k{{\rm soc}(G)}}(k\Omega)$ is always a symmetric algebra for arbitrary field $k$. Invoking Lemma 3.3, we deduce that ${\rm End}_{kG}(k\Omega)$ is always a symmetric algebra when $G$ is an almost simple primitive permutation group with prime power degree.

\medskip

To end up, we mention the following result as a consequence of the methods used in this paper. Recall that a block for a transitive permutation group $G$ on $\Omega$ is a subset $\Delta\subset\Omega$ such that either $\Delta^g=\Delta$ or $\Delta^g\cap\Delta=\emptyset$ for all $g\in G$.

\crllr Let $k$ be a field with characteristic $p$ and $G$ a transitive permutation group on $\Omega$ of degree $q^e$ such that $e\leq 2$ and $q\neq p$ is a prime, then ${\rm End}_{kG}(k\Omega)$ is a symmetric algebra provided a Sylow $q$-subgroup of $G$ is not normal when $q>p$.
\ecrllr
\medskip
\proof
First we suppose that $q<p$. Then it suffices to consider the case $e=2$ by Maschke's theorem. If $G$ is not primitive, then by [11, Theorem 16.1], $G$ contains an intransitive normal subgroup $N$ inducing $q$ blocks such that $G/N$ acts faithfully on these blocks. In particular, the assumption $q<p$ yields that $p\nmid |G/N|$ . Denote by $H$ a point stabilizer for some $\alpha\in\Omega$. Write $K=N\cap H$ and denote $$C={\rm Core}_N(K):=\cap_{n\in N} K^n.$$ Since $N/C$ acts faithfully on the block containing $\alpha$, we deduce that $p\nmid |N/C|$ and hence $p\nmid |H:C|$. Thus $k_H|k^H_C$, which means that $k^G_H$ is a direct summand of $k^G_C=(k^H_C)^G_H$. That $p\nmid |N/C||G/N|$ implies that $k^G_C=(k^N_C)^G_N$ is a semisimple $kG$-module, forcing $k^G_H$ to be a semisimple module too. Therefore $${\rm End}_{kG}(k\Omega)\cong {\rm End}_{kG}(k^G_H)$$ is a semisimple algebra, we are done in this case. Thus we may assume that $G$ is primitive. If $G$ is $2$-transitive, then we are in Situation 1 and we are done. Now we are left with the case $G$ is uniprimitive, ie., a primitive group which is not $2$-transitive. By [11, Theorem 16.3], a Sylow $q$-subgroup of $G$ is of order $q^2$ which is clearly an abelian regular subgroup. Invoking Lemma 3.2, we obtain the conclusion.

Next we suppose that $q>p$. Since a Sylow $q$-subgroup of $G$ is not normal, the result in [3, Theorem 3] implies that either $G$ is $2$-transitive on $\Omega$ or $G$ has a regular subgroup of order $q^e$. If $G$ is $2$-transitive, then analysis in Situation 1 yields the conclusion. If $G$ has a regular subgroup of order $q^e$, then the same reasoning as in the proof of Lemma 3.2 can be applied, as required.
\eprf

\end{document}